\documentclass[a4,12pt]{article}
\usepackage{amscd}
\usepackage{amsfonts}
\usepackage{amssymb}

\def\SU{\mathop{\mathrm {SU}}\nolimits}
\def\SL{\mathop{\mathrm {SL}}\nolimits}
\def\Mat{\mathop{\mathrm {Mat}}\nolimits}
\def\K{\mathop{\mathrm {K}}\nolimits}
\def\H{\mathop{\mathrm {H}}\nolimits}
\def\HP{\mathop{\mathrm {HP}}\nolimits}
\def\SO{\mathop{\mathrm {SO}}\nolimits}
\def\su{\mathop{\mathfrak {su}}\nolimits}
\def\Lie{\mathop{\mathrm {Lie}}\nolimits}
\def\End{\mathop{\mathrm {End}}\nolimits}
\def\Mat{\mathop{\mathrm {Mat}}\nolimits}

\newtheorem{thm}{Theorem}[section]

\begin{document}
\title{The noncommutative Chern-Connes character of the locally compact quantum normalizer of $\SU(1,1)$ in $\SL(2,\mathbb C)$}
\author{Do Ngoc Diep}
\maketitle
\begin{abstract}
We observe that the von Neumann (for short, W*-)envelope of the quantum algebra of functions on the normalizer of the group $\SU(1,1)\cong \SL(2,\mathbb R)$ in $\SL(2,\mathbb C)$ via deformation quantization contains the von Neumann algebraic quantum normalizer of $\SU(1,1)$ in the frame work of Waronowicz-Korogodsky. We then use the technique of reduction to the maximal subgroup to compute the K-theory, the periodic cyclic homology and the corresponding Chern-Connes character. 
\end{abstract}

\section*{Introduction}
It was remarked \cite{koelink_kustermans} that among the short list of very few well-studied locally compact quantum groups: quantum E(2), quantum ``ax+b'', quantum ``az+b'', quantum ${}_q\SU(1,1)$, the quantum group ${}_q\SU(1,1)$ plays an important role. The representation theory of ${}_q\SU(1,1)$ was well-treated and fully described by I. M. Burban and A. U. Klimyk \cite{burban_klimyk}, see also \cite{koelink_stokman}. 

In our previous works \cite{diepkukutho1}, \cite{diepkukutho2} and \cite{diepkuku} we developed a method of computation of the K-groups, periodic cyclic homology groups and the corresponding noncommutative Chern-Connes characters as homomorphisms between the two theories. Our method is based on the following ingredients:
\begin{enumerate}
\item Reduce the computation to the smooth case, the quantized algebras of smooth functions on Lie groups.
\item Reduce the computation to the case of maximal compact subgroups.
\item For the compact quantum groups, use the method developed in \cite{diepkukutho1}, \cite{diepkukutho2}.
\end{enumerate}
This method was applied in \cite{diepkuku} for the quantized algebras of functions on coadjoint orbits of the Lie groups like ``ax+b'', ``az+b'' and $\SL(2,\mathbb R)$ which are obtained from the deformation quantization of the algebras of functions on coadjoint orbits. 

In this paper we compute the K-theory groups, the periodic cyclic homology groups and the Chern-Connes character between them for the quantum groups ${}_q\widetilde{\SU(1,1)}$. We show that the K-theory and cyclic theory in this case are isomorphic to the corresponding ones for the torus $\mathbb T= \mathbb S^1$ and that the noncommutative Chern-Connes character is equivalent to the ordinary (cohomological) Chern character in the ordinary case of torus $\mathbb T=\mathbb S^1$.

In order to make clear the ideas and situation, we draw in the section 1  a corollary of the method of computation of K-theory, cyclic theory and noncommutative Connes-Chern character for the group C*-algebra $C^*(\SU(1,1))$. 
In Section 2 we prove that the deformation-quantized algebra   $\widetilde{\SU(1,1)}$ of smooth function with compact support can be included in the von Neumann algebraic quantum group ${W}^*_q(\widetilde{\SU(1,1)})$ as some dense subalgebra. Section 3  and Section 4 are devoted to reduction of the computation to the compact Lie group case and then reduction to the maximal compact subgroups.

\section{Noncommutative Chern-Connes character of the group C*-algebra $C^*(\SU(1,1))$}
Let us first recall that the normalizer of $\SU(1,1)$ in $\SL(2,\mathbb C)$ is the subgroup consisting of $2 \times 2$ complex matrices $X \in \Mat_2(\mathbb C)$ such that $X^*UX = \pm U$, where $U= \left[\begin{array}{cc} 1 & 0 \\ 0 & -1 \end{array}\right]$. 

The following result is an easy consequence of the main theorem of \cite{diepkuku} and \cite{diepkukutho1}.
\begin{thm}
The K-theory an the cyclic theory for $C^*(\SU(1,1))$ and the corresponding W-equivariant theories for $C(\mathbb T)$ are isomorphic, i.e.
$$\K_*(C^*(\SU(1,1)) \cong \K_*^W(C(\mathbb T)) \cong \K^{*,W}(\mathbb T),$$
$$\HP_*(C^*(\SU(1,1)) \cong \HP_*(C(\mathbb T)) \cong \H^{*,W}_{DR}(\mathbb T),$$
where $H^*_{DR}$  denote the $\mathbb Z_2$-graded de Rham cohomology of $\mathbb T$, and the noncommutative Chern-Connes character 
$$ch: \K_*(C^*(\SU(1,1))) \to \HP_*(C^*(\SU(1,1)))$$ becomes the ordinary Chern character of the torus and is therefore an isomorphism. 
\end{thm}
{\sc Proof.}
Following the main theorem of V. Nistor \cite{nistor} and our main theorem of \cite{diepkuku}, the K-theory and the cyclic theory of $C^*(\SU(1,1))$ are isomorphic with the same theories for the C*-algebra $C^*(\SO(2))$ of the maximal compact subgroup $\SO(2)$, which is also the maximal compact torus inside this maximal compact subgroup itself.  The K-theory and the cyclic theory for $\SO(2) \approx \mathbb S^1$ are isomorphic with the corresponding $W$-equivariant cohomological K-theory and $\mathbb Z_2$-graded de Rham theory. The noncommutative Chern-Connes character is equivalent to the classical Chern character of torus $\mathbb S^1$, i.e. we have a commutative diagram with vertical isomorphisms and the bottom isomorphism
$$\CD \K_*(C^*(\SU(1,1))) @>ch>> \HP_*(C^*(\SU(1,1)))\\
@V\cong VV @VV \cong V \\
K^*(\mathbb T) @>ch>\cong > H^*_{DR}(\mathbb T) \endCD$$
The Chern-Connes character is therefore an isomorphism. \hfill$\Box$

\section{The deformation-quantized algebra of functions as a dense subalgebra}
\begin{thm}\label{thm2.1}
The von Neumann envelop of the quantized algebra  $C^\infty_c(\widetilde{\SU(1,1)})$ of smooth functions via deformation quantization of the Poisson structure of the Lie algebra $\su(1,1)$ contains, and therefore is isomorphic to,  the von Neumann algebraic quantum algebra $W^*_q((\widetilde{\SU(1,1)})$ with generators $\alpha, \beta, \gamma, \delta$ subject to the relations
$$\alpha \beta = q \beta \alpha, \alpha \gamma = q \gamma \alpha, \beta \delta = q \delta\beta, \gamma\delta = q \delta\gamma,$$
$$\beta\gamma = \gamma\beta, \alpha\delta - q \beta\gamma = \delta\alpha - q^{-1}\beta\gamma =1,$$ where $1$ denotes the unit in $W^*_q((\widetilde{\SU(1,1)})$ and $0<q<1$,
as a dense subalgebra.
\end{thm}
{\sc Proof.}
Let us denote $x_{ij}$ the matrix coefficients of the standard representation of $\SU(1,1)$ in $\mathbb C^2$: Let $X_+, X_-, K_\pm$ be the natural basis of $\mathfrak g = \su(1,1)$ and $\rho$ the standard representation of $U_h(\su(1,1))$ given by
$$X_+ \mapsto \left[\begin{array}{cc} 0 & 1\\ 0 & 0 \end{array}\right]$$
$$X_- \mapsto \left[\begin{array}{cc} 0 & 0\\ 1 & 0 \end{array}\right]$$
$$K_\pm \mapsto \left[\begin{array}{cc} q^{\pm 1} & 0 \\ 0 & q^{\mp 1}\end{array}\right]$$ 
Let $\Delta$ be the product of $U_h(\su(1,1))$, e. i. 
$$\Delta(X_+) = X_+ \otimes K_+ + 1 \otimes X_+,$$
$$\Delta(X_-) = X_- \otimes 1 + K_- \otimes X_-,$$
$$\Delta(K_\pm) = K_\pm \otimes K_\pm.$$

Then the quantized universal enveloping algebra $U_h(\su(1,1))$ is isomorphic as  $\mathbb C[[h]]$-modules to $U(\su(1,1))$. The convolution product can be defined by 
$$f \star g(x) := f \otimes g (\Delta(x)),$$ and therefore
$$g \star f (x)= g \otimes f(\Delta(x)) = f \otimes g (\Delta^{op}(x)) = f \otimes g(R \Delta (R^{-1})(x))$$ where by definition 
$$R := e^{\frac{h}{8}{H \otimes H}} \sum\frac{(1-q^{-2})q^{1/2n(n-1)}}{[n]_q }(K_+X_+ \otimes K_-X_-)^n,$$ and $[n]_q := \frac{q^n -q^{-n}}{q - q^{-1}}$. It is well-known, see e.g. \cite{shnider_sternberg}
that $$x_{ij} \star x_{kl} = x_{kl} \otimes x_{ij}(R\Delta R^{-1}) = R^{ki}_{pq}x^p_i \star x^q_s (R^{-1})^{rs}_{lj}$$ where $R$ is the tensor 
$$R= q^{1/2}\left[\begin{array}{cccc}
q & 0 & 0 & 0 \\
0 & 1 & q & 0 \\
0 & q^{-1} & 1 & 0\\
0 & 0 & 0 & q \end{array}\right] $$ in $\End(V\otimes V)$, $V = \mathbb C^2$.
If we put $\alpha = x^1_1$, $x^1_2 = \beta$, $x^2_1 = \gamma$ and $x^2_2 = \delta$ we have relations
$$\alpha \star \beta = q \beta \star \alpha, \alpha \star \gamma = q \gamma \star \alpha, \beta \star \delta = q \delta \star \beta, \gamma \star \delta = q \gamma \star \delta, $$
$$ \beta \star \gamma = \gamma \star \beta, \alpha \star \delta - q \beta \star \delta = \delta \star \alpha - q^{-1}\beta \star \gamma = 1.$$
This means that the matrix coefficients $x_{ij}$ of the standard representation $\rho$ generated the space of all polynomial functions.
 
It was shown, see e.g. (\cite{shnider_sternberg}, Propositions 11.9.1, 11.9.2) that the bialgebra map from $W^*_q(\widetilde{\SU(1,1)}) $ to $U_h(\su(1,1))$ is injective on generators and that the algebra
$W^*_q(\widetilde{\SU(1,1)})$ is the quantized Poisson-Lie algebra associated to the Lie bialgebra which is the classical limit of $U_h(\su(1,1))$.

 Elements of the universal enveloping algebra $U(\su(1,1))$ can be considered also as polynomial functions over the dual space to the Lie algebra $\su(1,1)$. Let us recall the Possoin structure on the dual space of the Lie algebra $\su(1,1) = \Lie \SU(1,1)$ : for all $f,g\in C^\infty(\su(1,1))$
their Possoin bracket is $$\{f,g\}(F) := \langle F,[df,dg]\rangle,$$  for all $F \in \mathfrak g^* = \su(1,1)^*$, where $df, dg \in Hom(\mathfrak g^*, \mathbb R) \cong \mathfrak g$. 
To this structure associates a star product $*$.

Apply this construction for $f= x_{ij}, g= x_{kl}$ 

We have $$\{f,g\} = \frac{1}{h}(x\star y - y \star x) (mod\; h^2)$$ Let us consider now the ordinary star product of functions on the Lie algebra $\widetilde{\su(1,1)}$. Denote again the standard representation by $\rho : \widetilde{\SU(1,1)} \to \Mat_2(\mathbb C)$ and the matrix coefficients satisfy the orthogonal rations,  we have therefore the relation
$\rho(f * g) = \rho(f)\rho(g)$. Because the standard representation is faithful, we can deduce that 
$$\Delta(x_{ij}) = \sum_k x_{ik} \otimes x_{kj}$$ and
$$f\star g = (f \otimes g)\circ \Delta = f * g.$$
Thus, the generators of von Neumann algebraic $W^*_q(\SU(1,1))$ are in a bijection with the functions $x_{ij}$ and the product structure are agreed. So the von Neumann envelop of $C^\infty(\SU(1,1))$ contains and therefore isomorphic to the $W^*_q(\SU(1,1))$.   
The theorem is therefore proven.
\hfill$\Box$

\section{Restriction to a maximal compact subgroup}
Now we use the technique of reduction to maximal compact subgroups developed in \cite{nistor} and \cite{diepkukutho1}. 
\begin{thm}\label{thm3.1}
The K-theory and the cyclic theory for $W^*_q(\widetilde{\SU(1,1)})$ are isomorphic to the corresponding theories for $C^\infty_c(\SO(2))$, i.e.
$$\K_*(W^*_q(\widetilde{\SU(1,1)})) \cong \K_*(C^\infty_c(\widetilde{\SU(1,1)}))\cong \K_* C^\infty_c(\widetilde{\SO(2)}),$$  
$$\HP_*(W^*_q(\widetilde{\SU(1,1)})) \cong \HP_*(C^\infty_c(\widetilde{\SU(1,1)}))\cong \HP_* C^\infty_c(\widetilde{\SO(2)}).$$  
\end{thm}
{\sc Proof.}
From the general cyclic theory of A. Connes, one reduces
$$\K_*(W^*_q(\widetilde{\SU(1,1)})) \cong \K_*(C^\infty_c(\widetilde{\SU(1,1)}))$$ and 
$$\HP_*(W^*_q(\widetilde{\SU(1,1)})) \cong \HP_*(C^\infty_c(\widetilde{\SU(1,1)})).$$ From the results of \cite{nistor} and \cite{diepkuku} we have the second isomorphisms
$$\K_*(C^\infty_c(\widetilde{\SU(1,1)}))\cong \K_* C^\infty_c(\widetilde{\SO(2)}) $$ and
$$\HP_*(C^\infty_c(\widetilde{\SU(1,1)}))\cong \HP_* C^\infty_c(\widetilde{\SO(2)}).$$
\hfill$\Box$

\section{Noncommutative Chern-Connes character}
Let us finally draw back the corresponding computation results and in particular the noncommutative Chern-Connes character.
\begin{thm}
Let us denote $\tilde{\mathbb T} = \mathbb T \rtimes \mathbb Z_2$ a fixed maximal compact subgroup in the normalizer $\widetilde{\SU(1,1)}$ of $\SU(1,1)$ in $\SL(2,\mathbb C)$. Then
$$ \K_*(W^*_q(\widetilde{\SU(1,1)})) \cong \K^{*,W}(\tilde{\mathbb T}),$$
$$ \HP_*(W^*_q(\widetilde{\SU(1,1)})) \cong \HP^{*,W}(\tilde{\mathbb T}),$$
and the noncommutative Chern-Connes character
$$ch : \K_*(W^*_q(\widetilde{\SU(1,1)})) \to \HP_*(W^*_q(\widetilde{\SU(1,1)}))$$ is an isomorphism, where $W$ denotes the Weyl group corresponding to the torus.
\end{thm}
{\sc Proof.}
The proof is a combination of Theorems \ref{thm2.1} and \ref{thm3.1} and the main result of \cite{diepkukutho2}
 The K-theory and the cyclic theory for $\SO(2) \approx \mathbb S^1$ are isomorphic with the corresponding cohomological K-theory and $\mathbb Z_2$-graded de Rham theory. The Chern-Connes character is equivalent to the classical Chern character of torus $\mathbb S^1$, i.e. we have a commutative diagram with vertical isomorphisms and the bottom isomorphism
$$\CD \K_*(W^*_q(\widetilde{\SU(1,1)})) @>ch>> \HP_*(W^*_q(\widetilde{\SU(1,1)}))\\
@V\cong VV @VV \cong V \\
K^{*,W}(\widetilde{\mathbb T}) @>ch>\cong > H^{*,W}_{DR}(\widetilde{\mathbb T}) \endCD$$
The Chern-Connes character is therefore an isomorphism.
\hfill$\Box$

\section{Acknowledgments}
The main part of the paper was realized during a stay of the author in Abdus Salam ICTP. The author would like to express sincere thanks to ICTP, and in particular Professor Le Dung Trang and Professor Aderemi O. Kuku, for invitation and for the provided excellent conditions of work. The deep thanks are addressed to professor E. Koelink for the useful remarks, and especially for the reference [KS].

{\noindent\sc Institute of Mathematics, National Center for Science and Technology of Vietnam, 18 Hoang Quoc Viet Road, Cau Giay District, 10307 Hanoi, Vietnam}\\
{\tt e-Mail: dndiep@math.ac.vn}
\end{document}